\documentclass[10pt,draft,twoside,leqno,a4paper]{article}
\usepackage{amsmath,amsthm,amsfonts,amssymb,verbatim}
\newtheorem{Theorem}{Theorem}[section]
\newtheorem{Lemma}{Lemma}[section]

\newtheorem{Proposition}{Proposition}[section]

\newcommand{\R}{\mathbb R}
\newcommand{\eps}{\varepsilon}

\newcommand{\dx}{\,{\mathrm d}x}
\newcommand{\dy}{\,{\mathrm d}y}
\newcommand{\weak}{\rightharpoondown}
\newcommand{\bra}{\langle}
\newcommand{\ket}{\rangle}
\newcommand{\const}{A}
\newcommand{\us}{{\widetilde u}}

\newcommand{\sing}{\sigma}
\newcommand{\Ge}{G_\varepsilon}
\newcommand{\ue}{ u_\varepsilon}
\newcommand{\uue}{\underline{u}_\varepsilon}

\newcommand{\tue}{\widetilde u_\varepsilon}
\newcommand{\tuo}{\widetilde u_0}
\newcommand{\Ae}{\mathcal A_\varepsilon}
\newcommand{\Ao}{\mathcal A_0}

\newcommand{\Ie}{I_\varepsilon}

\newcommand{\uo}{u_0}
\newcommand{\uuo}{\underline{u}_0}
\newcommand{\Io}{I_0}

\newcommand{\Ue}{U_\varepsilon}
\newcommand{\Fe}{F_\varepsilon}
\newcommand{\Fo}{F_0}
\newcommand{\we}{w_\varepsilon}
\newcommand{\uj}{{u_j}}
\newcommand{\ujp}{u_j'}
\newcommand{\cj}{c_j}

\begin{document}
\title{Multiplicity for a nonlinear fourth order
elliptic equation in Maxwell-Chern-Simons vortex theory}
\author{Tonia Ricciardi\thanks{
Partially supported by MIUR, National Project ``Variational Methods and Nonlinear
Differential Equations"}\\
{\small Dipartimento di Matematica e Applicazioni}\\
{\small Universit\`a di Napoli Federico II}\\
{\small Via Cintia}\\
{\small 80126 Naples, Italy}\\
{\small fax: +39 081 675665}\\
{\small e-mail: tonia.ricciardi@unina.it}
}
\maketitle
\begin{abstract}
We prove the existence of at least two solutions for a fourth
order equation, which includes the vortex equations for the $U(1)$
and $CP(1)$ self-dual Maxwell-Chern-Simons models as special
cases. Our method is variational, and it relies on an ``asymptotic
maximum principle" property for a special class of supersolutions
to this fourth order equation.
\end{abstract}
\begin{description}
\item{\textsc{Key Words:}} nonlinear fourth-order elliptic
equation, nonlinear elliptic system, maximum principle,
Chern-Simons vortex theory \item{\textsc{MSC 2000 Subject
Classification:}} 35J60
\end{description}
\setcounter{section}{-1}
\section{Introduction}
\label{sec:introduction} Vortex solutions for self-dual
Maxwell-Chern-Simons  models may be generally reduced to systems
of two nonlinear elliptic equations of the second order, defined
on two-dimensional Riemannian manifolds. See \cite{KLL,LLM} and
the monographs \cite{D,JT,Y}. These systems are also equivalent to
scalar nonlinear elliptic equations of the fourth order. The
existence of {\em multiple} solutions for such fourth-order
equations, in the case of compact manifolds, is the main question
addressed in this note.
\par
We were motivated to consider this problem by our previous joint
work with Tarantello \cite{RT} concerning the $U(1)$
Maxwell-Chern-Simons model introduced in \cite{LLM}, and by the
results of Chae and Nam \cite{CN} concerning the $CP(1)$
Maxwell-Chern-Simons model introduced in \cite{KLL}. By
variational methods, it is shown in \cite{RT} that the $U(1)$
system admits in general at least two distinct vortex solutions.
On the other hand, the method employed in \cite{CN} allows the
authors to obtain only one vortex solution for the $CP(1)$ system.
\par
Our main result (see Theorem~\ref{thm:main} below) will yield
multiple solutions for a {\em general} system containing the
$U(1)$ system and the $CP(1)$ system as special cases. As in the
$U(1)$ case, we shall reduce the system to an elliptic fourth
order equation admitting a convenient variational structure. We
shall obtain two solutions corresponding to a local minimum and a
mountain pass. Due to our abstract formulation, we cannot use the
{\em ad hoc} methods employed in \cite{RT}, based on minimization
with integral constraints. Instead, we shall exploit an
``asymptotic maximum principle" property for a special class of
supersolutions of this fourth order equation, which we believe is
of interest of its own.
\par
More precisely, we denote by $M$ a compact Riemannian 2-manifold,
and we fix $n>0$ points $p_1,\ldots,p_n\in M$ (we already showed
in \cite{R2} that the case $n=0$ admits only the trivial solution
$(e^\us,v)=(f^{-1}(s),s)$). We consider (distributional) solutions
$(\us,v)$ for the system:
\begin{align}
\label{genmcsa'}
&-\Delta\us=\eps^{-1}\lambda(v-f(e^\us))-4\pi\sum_{j=1}^n\delta_{p_j}
&&\text{on}\ M\\
\label{genmcsb'} &-\Delta v=\eps^{-1}\left[\lambda
f'(e^\us)e^\us(s-v)-\eps^{-1}(v-f(e^\us))\right] &&\text{on}\ M.
\end{align}
Here $s\in\R$, $\eps>0$, $\lambda>0$ are constants,
$\delta_{p_j}$, $j=1,\ldots,n$ is the Dirac measure centered at
$p_j$, and $f:[0,+\infty)\to\R$ is smooth and {\em strictly
increasing}, i.e., $f'(t)>0$ for all $t\in[0,+\infty)$. Some
further technical assumptions on $f$ will be made below. System
\eqref{genmcsa'}--\eqref{genmcsb'} was introduced in \cite{R2}. It
contains the system for $U(1)$ Maxwell-Chern-Simons vortices
introduced in \cite{LLM} and analyzed in \cite{CK,R1,RT}, as well
as the system for $CP(1)$ Maxwell-Chern-Simons vortices (in the
``single-signed case") introduced in \cite{KLL} and analyzed in
\cite{CN} as special cases. Indeed, the $U(1)$ system \cite{LLM}
is given by:
\begin{align*}
&\Delta\us=2q^2e^\us-2\mu N+4\pi\sum_{j=1}^n\delta_{p_j}
&&\text{on}\ M\\
&\Delta
N=(\mu^2+2q^2e^\us)N-q^2(\mu+\frac{2q^2}{\mu})e^\us&&\text{on}\ M.
\end{align*}
Setting $\lambda=2q^2/\mu$, $\eps=1/\mu$, $v:=\mu/q^2\,N$, the
above system takes the form:
\begin{align}
\label{u1a'}
-\Delta\us=&\eps^{-1}\lambda(v-e^\us)-4\pi\sum_{j=1}^n\delta_{p_j}
&&\text{on}\ M\\
\label{u1b'} -\Delta v=&\eps^{-1}\{\lambda
e^\us(1-v)-\eps^{-1}(v-e^\us)\} &&\text{on}\ M,
\end{align}
which corresponds to \eqref{genmcsa'}--\eqref{genmcsb'} with
$f(t)=t$ and $s=1$. On the other hand, the $CP(1)$ system
\cite{KLL} is given by:
\begin{align*}
&\Delta\us=2q\left(-N+S-\frac{1-e^\us}{1+e^\us}\right)+4\pi\sum_{j=1}^n\delta_{p_j}
&&\text{on}\ M\\
 &\Delta
N=-\kappa^2q^2\left(-N+S-\frac{1-e^\us}{1+e^\us}\right)+q\frac{4e^\us}{(1+e^\us)^2}N
&&\text{on}\ M.
\end{align*}
Setting $v=N-S$, $s=-S$, $\lambda=2/\kappa$, $\eps=1/(\kappa q)$,
the system above takes the form:
\begin{align}
\label{cp1a'}
-\Delta\us=&\eps^{-1}\lambda\left(v-\frac{e^\us-1}{e^\us+1}\right)-4\pi\sum_{j=1}^n\delta_{p_j}
&&\text{on}\ M\\
\label{cp1b'} -\Delta
v=&\eps^{-1}\Big\{\lambda\frac{2e^\us}{(1+e^\us)^2}(s-v)-\eps^{-1}
\left(v-\frac{e^\us-1}{e^\us+1}\right)\Big\} &&\text{on}\ M,
\end{align}
which corresponds to \eqref{genmcsa'}--\eqref{genmcsb'} with
$f(t)=(t-1)/(t+1)$.
\par
In turn, system \eqref{genmcsa'}--\eqref{genmcsb'} is equivalent
to the following {\em fourth order} equation (see
Section~\ref{sec:variational} for the details):
\begin{align}
\nonumber \eps^2\Delta^2 u-&\Delta u=-\eps\lambda
[f''(e^{\sing+u})e^{\sing+u}+f'(e^{\sing+u})]e^{\sing+u}|\nabla(\sing+u)|^2\\
\label{fourthorderintro} &+2\eps\lambda\Delta f(e^{\sing+u})
+\lambda^2f'(e^{\sing+u})e^{\sing+u}(s-f(e^{\sing+u}))-\frac{4\pi
n}{|M|}\qquad\text{on}\ M,
\end{align}
where $\sing$ is the Green function uniquely defined by
$-\Delta\sing=4\pi(n/|M|-\sum_{j=1}^n\delta_{p_j})$ on $M$ and
$\int_M\sing=0$, with $|M|$ the volume of $M$.
\par
We make the following
\par
\medskip\noindent \textbf{Assumptions on $f$:}
\begin{enumerate}
\item[($f$0)] $f:[0,+\infty)$ is smooth and $f'(t)>0$ for all
$t\in[0,+\infty)$ \item[($f$1)] $f(0)<s<\sup_{t>0}f(t)$
\item[($f$2)] $f$, $f'$, $f''$ have at most polynomial growth
\item[($f$3)] $f$ satisfies one of the following conditions:
\begin{enumerate}
\item[(a)]$f''(t)t+f'(t)\ge0$ and
$\sup_{t>0}|f(t)|/[f'(t)t]<+\infty$
\item[(b)]$\sup_{t>0}f'(t)t(|\log t|+|f(t)|)<+\infty$.
\end{enumerate}
\end{enumerate}
This aim of this note is to establish the following result for
\eqref{fourthorderintro}:
\begin{Theorem}
\label{thm:main} Suppose $f$ satisfies assumptions ($f$0), ($f$1),
($f$2) and ($f$3). Then there exists $\lambda_0>0$ with the
property that for every $\lambda\ge\lambda_0$ there exists
$\eps_\lambda>0$ such that the fourth order equation
\eqref{fourthorderintro} admits at least {\em two} solutions for
all $0<\eps<\eps_\lambda$.
\end{Theorem}
We note that assumption ($f$3)--(a) allows $f(t)=t^\alpha$, for
every $\alpha>0$, and therefore it includes the $U(1)$ case
$f(t)=t$. On the other hand, assumption ($f$3)--(b) is satisfied
by the $CP(1)$ case $f(t)=(t-1)/(t+1)$. It follows that the
existence result stated in Theorem~\ref{thm:main} includes indeed
the $U(1)$ system and the $CP(1)$ system as special cases, as well
as all power growths for $f$.
\par
 As already mentioned, we shall prove
Theorem~\ref{thm:main} variationally. Indeed, in
Section~\ref{sec:variational} we show  that solutions to
\eqref{fourthorderintro} correspond to critical points for the
functional
\begin{align*}
\Ie(u)=\frac{\eps^2}{2}&\int(\Delta u)^2+\frac{1}{2}\int|\nabla u|^2\\
+&\eps\lambda\int f'(e^{\sing+u})e^{\sing+u}|\nabla(\sing+u)|^2
+\frac{\lambda^2}{2}\int(f(e^{\sing+u})-s)^2+\frac{4\pi
n}{|M|}\int u,
\end{align*}
defined on the Sobolev space $H^2(M)$ (we choose to emphasize the
dependence on $\eps$ only, since $\lambda$ will be fixed in
Section~\ref{sec:subsol}). In Section~\ref{sec:subsol} we show
that for $\eps$ and $\lambda$ as in Theorem~\ref{thm:main}, the
functional $\Ie$ admits a local minimum $\uue$. This subsolution
will be employed in Section~\ref{sec:locmin} to obtain a critical
point $\ue$ satisfying $\Ie(\ue)=\min\{\Ie(u)/u\in
H^2(M),u\ge\uue\}$, in the spirit of some results of Brezis and
Nirenberg \cite{BN} concerning second order equations. Since our
equation is of the fourth order, and thus the Hopf maximum
principle is not directly applicable, the main technical
difficulty will be to show that $\ue>\uue$, pointwise on $M$.
Nevertheless, by exploiting the decomposition
$\eps^2\Delta^2-\Delta=(-\eps^2\Delta+1)(-\Delta)$, we shall
derive a kind of ``strong maximum principle" property for a
special class of supersolutions of \eqref{fourthorderintro}, for
small values of $\eps$. Finally, in Section~\ref{sec:PS} we show
that under the assumptions of Theorem~\ref{thm:main}, the
functional $\Ie$ satisfies the Palais-Smale condition. Therefore,
the existence of a second critical point of the ``mountain pass"
type will follow by the Ambrosetti and Rabinowitz theorem
\cite{AR}.
\par
We showed in our previous note \cite{R2} that under assumptions
($f$0) and ($f$1), solutions for system
\eqref{genmcsa}--\eqref{genmcsb} (equivalently solutions to
\eqref{fourthorderintro}) tend to a solution for
\begin{equation}
\label{gencs} -\Delta
u=\lambda^2f'(e^{\sing+u})e^{\sing+u}(s-f(e^{\sing+u}))-\frac{4\pi
n}{|M|}\qquad\text{on}\ M,
\end{equation}
as $\eps\to0$, in any relevant topology. Equation \eqref{gencs} is
a generalization of the $U(1)$ ``pure" Chern-Simons equation
derived in \cite{HKP,JW} and thoroughly analyzed in
\cite{CY,T,NT,DJLW} in the case of compact manifolds (see
references therein for the non-compact case), and of the $CP(1)$
``pure" Chern-Simons equation derived in \cite{KLL1} and analyzed
in \cite{CN1}. We note that solutions for \eqref{gencs} correspond
to critical points for the functional $\Io$ defined for $u\in
H^1(M)$ by
\[
\Io(u)=\frac{1}{2}\int|\nabla
u|^2+\frac{\lambda^2}{2}\int(f(e^{\sing+u})-s)^2 +\frac{4\pi
n}{|M|}\int u.
\]
A multiplicity result may be obtained for \eqref{gencs}
variationally, which implies that the multiplicity results for the
particular ``pure'' Chern-Simons equations obtained in
\cite{CN1,T} are in fact a {\em general property} of equations of
the form \eqref{gencs}. Since the proof may be obtained from the
proof of Theorem~\ref{thm:main} by setting $\eps=0$, we omit the
details.
\par
{\em Notation.} Henceforth, unless otherwise specified, all
equations are defined on $M$. All integrals are taken over $M$
with respect to the Lebesgue measure. All functional spaces are
defined on $M$ in the usual way. In particular, we denote by
$L^p$, $1\le p\le+\infty$, the Lebesgue spaces and by $H^k$, $1\le
k\le +\infty$, the Sobolev spaces. We denote by $C>0$ a general
constant, independent of certain parameters that will be specified
in the sequel, and whose actual value may vary from line to line.
\section{Variational setting}
\label{sec:variational} Our aim in this section is to provide a
suitable variational formulation for the generalized
Maxwell-Chern-Simons system \eqref{genmcsa'}--\eqref{genmcsb'}, by
reducing it to the fourth order equation \eqref{fourthorderintro}.
We note that the resulting formulation is new even for the special
case $f(t)=(t-1)/(t+1)$, corresponding to the $CP(1)$ model.
\par
In order to work in Sobolev spaces, it is standard (see
\cite{JT,Y}) to subtract from $\us$ its ``singular part", which we
denote by $\sing$. Namely, we denote by $\sing$ the unique
solution for the problem
\begin{align*}
&-\Delta\sing=A-4\pi\sum_{j=1}^n\delta_{p_j}\\
&\int\sing=0,
\end{align*}
where $A=4\pi n/|M|>0$. Since $\sing(x)\approx\log|x-p_j|^2$ near
$p_j$, $j=1,\ldots,n$, it follows that $e^\sing$ and
$e^\sing\nabla\sing=\nabla e^\sing$ are {\it smooth} on $M$.
Furthermore, $e^{\sing}|\nabla \sing|^2$ and $e^\sing\Delta\sing$
are also smooth. Indeed, it is easy to check that
\begin{align}
\label{singdelta}
&e^\sing\Delta\sing=-Ae^\sing\\
\label{singnabla} &e^{\sing}|\nabla \sing|^2=\Delta
e^{\sing}+\const e^\sing,
\end{align}
in the sense of distributions.
\par
Setting $\us=\sing+u$, we obtain from
\eqref{genmcsa'}--\eqref{genmcsb'} the equivalent system for
$(u,v)\in H^1(M)\times H^1(M)$:
\begin{align}
\label{genmcsa}
&-\Delta u=\eps^{-1}\lambda\left(v-f(e^{\sing+u})\right)-\const\\
\label{genmcsb} &-\Delta v= \eps^{-1}\left[\lambda
f'(e^{\sing+u})e^{\sing+u}(s-v)
-\eps^{-1}\left(v-f(e^{\sing+u})\right)\right].
\end{align}
System \eqref{genmcsa}--\eqref{genmcsb} is equivalent to a fourth
order equation:
\begin{Lemma}
\label{lem:fourthorder} $(u,v)\in H^1\times H^1$ is a weak
solution for \eqref{genmcsa}--\eqref{genmcsb} if and only $u\in
H^2$ is a weak solution for the fourth order equation
\begin{align}
\nonumber \eps^2\Delta^2 u-&\Delta u=-\eps\lambda
[f''(e^{\sing+u})e^{\sing+u}+f'(e^{\sing+u})]e^{\sing+u}|\nabla(\sing+u)|^2\\
\label{fourthorder} &+2\eps\lambda\Delta f(e^{\sing+u})
+\lambda^2f'(e^{\sing+u})e^{\sing+u}(s-f(e^{\sing+u}))-\const
\end{align}
and $v$ is defined by
\begin{equation}
\label{vdef} v=-\eps\lambda^{-1}\Delta
u+\eps\lambda^{-1}\const+f(e^{\sing+u}).
\end{equation}
\end{Lemma}
\begin{proof}
By elliptic regularity, weak solutions $(u,v)\in H^1\times H^1$
for \eqref{genmcsa}--\eqref{genmcsb} are smooth. Clearly,
\eqref{vdef} is equivalent to \eqref{genmcsa}. Inserting
\eqref{vdef} into \eqref{genmcsb}, we obtain
\begin{align*}
\eps\lambda^{-1}\Delta^2 u-\Delta f(e^{\sing+u})
=&\eps^{-1}\lambda
f'(e^{\sing+u})e^{\sing+u}(s+\eps\lambda^{-1}\Delta
u-\eps\lambda^{-1}A-f(e^{\sing+u}))\\
&-\eps^{-2}(-\eps\lambda^{-1}\Delta u+\eps\lambda^{-1}A).
\end{align*}
Equivalently, multiplying by $\eps\lambda$:
\begin{align}
\label{fourthprelim} \eps^2\Delta^2u-\Delta u=&\eps\lambda\Delta
f(e^{\sing+u})
+\eps\lambda f'(e^{\sing+u})e^{\sing+u}(\Delta u-A)\\
\nonumber
&+\lambda^2f'(e^{\sing+u})e^{\sing+u}(s-f(e^{\sing+u}))-A.
\end{align}
By \eqref{singdelta} we have:
\[
f'(e^{\sing+u})e^{\sing+u}(\Delta u-A)
=f'(e^{\sing+u})e^{\sing+u}\Delta(\sing+u).
\]
Furthermore, we have:
\begin{equation}
\label{deltafid} \Delta f(e^{\sing+u})
=\{f''(e^{\sing+u})e^{\sing+u}+f'(e^{\sing+u})\}e^{\sing+u}|\nabla(\sing+u)|^2
+f'(e^{\sing+u})e^{\sing+u}\Delta(\sing+u).
\end{equation}
Therefore,
\begin{align*}
\Delta f(e^{\sing+u})+&f'(e^{\sing+u})e^{\sing+u}(\Delta u-A)
=\Delta f(e^{\sing+u})+f'(e^{\sing+u})e^{\sing+u}\Delta(\sing+u)\\
=&2\Delta f(e^{\sing+u})-\{f''(e^{\sing+u})e^{\sing+u}
+f'(e^{\sing+u})\}e^{\sing+u}|\nabla(\sing+u)|^2,
\end{align*}
and \eqref{fourthprelim} reduces to \eqref{fourthorder}. The
converse follows similarly.
\end{proof}
Now we obtain a variational formulation for \eqref{fourthorder}:
\begin{Lemma}
\label{lem:variational} $u\in H^2$ is a weak solution for
\eqref{fourthorder} if and only if it is a critical point for the
functional $\Ie$ defined on $H^2$ by:
\begin{align*}
\Ie(u)=\frac{\eps^2}{2}&\int(\Delta u)^2+\frac{1}{2}\int|\nabla u|^2\\
+&\eps\lambda\int f'(e^{\sing+u})e^{\sing+u}|\nabla(\sing+u)|^2
+\frac{\lambda^2}{2}\int(f(e^{\sing+u})-s)^2+\const\int u.
\end{align*}
\end{Lemma}
\begin{proof}
First of all, we note that $\Ie$ is well-defined on $H^2$ and
smooth. Indeed, if $u\in H^2$, then by Sobolev embeddings we have
$|\nabla u|\in L^p$ for all $p\ge1$ and $u\in L^\infty$. We can
write, in view of \eqref{singnabla}
\begin{align}
\label{welldef}
\int f'(e^{\sing+u})&e^{\sing+u}|\nabla(\sing+u)|^2=\\
\nonumber =&\int
f'(e^{\sing+u})e^{\sing+u}(|\nabla\sing|^2+2\nabla\sing\cdot\nabla
u+|\nabla u|^2)\\
\nonumber =&\int f'(e^{\sing+u})e^u(\Delta e^\sing+Ae^\sing)+2\int
f'(e^{\sing+u})e^u\nabla e^\sing\cdot\nabla u\\
\nonumber &\qquad+\int f'(e^{\sing+u})e^{\sing+u}|\nabla u|^2,
\end{align}
and therefore the third term in $\Ie$ is well-defined on $H^2$.
Clearly, all other terms in $\Ie$ are well-defined, and thus $\Ie$
is well-defined on $H^2$. Smoothness of $\Ie$ is checked
similarly. We check that solutions in $H^2$ to \eqref{fourthorder}
correspond to critical points for $\Ie$. We compute, for any
$\phi\in H^2$:
\begin{align*}
\frac{d}{dt}\big|_{t=0}
\int f'(e^{\sing+u+t\phi})&e^{\sing+u+t\phi}|\nabla(\sing+u+t\phi)|^2\\
=&\int\left[f''(e^{\sing+u})e^{\sing+u}+f'(e^{\sing+u})\right]
e^{\sing+u}|\nabla(\sing+u)|^2\phi\\
&+2\int f'(e^{\sing+u})e^{\sing+u}\nabla(\sing+u)\cdot\nabla\phi,
\end{align*}
and therefore
\begin{align}
\label{gateaux}
\bra\Ie'(u),\phi\ket=
\eps^2&\int\Delta u\Delta\phi+\int\nabla u\cdot\nabla\phi\\
\nonumber
&+\eps\lambda\int\left[f''(e^{\sing+u})e^{\sing+u}+f'(e^{\sing+u})\right]
e^{\sing+u}|\nabla(\sing+u)|^2\phi\\
\nonumber
&+2\eps\lambda\int f'(e^{\sing+u})e^{\sing+u}\nabla(\sing+u)\cdot\nabla\phi\\
\nonumber
&+\lambda^2\int f'(e^{\sing+u})e^{\sing+u}(f(e^{\sing+u})-s)\phi
+\const\int\phi.
\end{align}
Since
\begin{align*}
\int f'(e^{\sing+u})e^{\sing+u}\nabla(\sing+u)\cdot\nabla\phi
=\int\nabla f(e^{\sing+u})\cdot\nabla\phi
=-\int\Delta f(e^{\sing+u})\phi,
\end{align*}
it follows that critical points of $\Ie$ correspond to solutions
to \eqref{fourthorder}, as asserted.
\end{proof}
We say that $u\in H^2$ is a subsolution (supersolution) for
\eqref{fourthorder} if $u$ satisfies $\bra\Ie'(u),\phi\ket\le0$
($\bra\Ie'(u),\phi\ket\ge0$), for all $\phi\in H^2$ such that
$\phi\ge0$ on $M$.
\section{Existence of a subsolution}
\label{sec:subsol}In this section we show that for suitable values
of $\lambda$ and $\eps$, the fourth order equation
\eqref{fourthorder} admits a subsolution. Namely, we show:
\begin{Proposition}
\label{prop:subsol} Suppose $f$ satisfies ($f$0) and ($f$1). Then
there exists $\lambda_0>0$ such that for every fixed
$\lambda\ge\lambda_0$ there exists $\eps_\lambda>0$ with the
property that for every $0<\eps<\eps_\lambda$ equation
\eqref{fourthorder} admits a subsolution $\uue$. Furthermore,
$\uue\to\uuo$ in $H^2$ as $\eps\to0$, with $\uuo$ a subsolution
for \eqref{gencs}.
\end{Proposition}
We begin by proving some properties of the Green function $\Ge$
for the operator $-\eps^2\Delta+1$:
\begin{Lemma}
\label{lem:green} Let $\Ge=\Ge(x,y)$ be the Green function defined
by
\[
(-\eps^2\Delta_x+1)\Ge(x,y)=\delta_y\qquad\text{on}\ M.
\]
Then
\begin{itemize}
\item[(i)] $\Ge>0$ on $M\times M$ and for every fixed $y\in M$ we
have $\Ge\weak\delta_y$ as $\eps\to0$, weakly in the sense of
measures; \item[(ii)] $\|\Ge\ast h\|_q\le\|h\|_q$ for all $1\le
q\le+\infty$; \item[(iii)] If $\Delta h\in L^q$ for some
$1<q<+\infty$, then $\|\Ge\ast h-h\|_q\le\eps^2\|\Delta h\|_q$.
\end{itemize}
\end{Lemma}
\begin{proof}
Proof of (i). Note that since $-\eps^2\Delta+1$ is coercive, $\Ge$
is well defined (e.g., by Stampacchia's duality argument
\cite{St}). By the maximum principle, $\Ge>0$ on $M\times M$.
Integrating over $M$ with respect to $x$, we have
$\int\Ge(x,y)\dx=\int|\Ge(x,y)|\dx=1$ and therefore there exists a
Radon measure $\mu$ such that $\Ge(\cdot,y)\weak\mu$ as
$\eps\to0$, weakly in the sense of measures. For $\varphi\in
C^\infty$ we compute:
\[
\varphi(y)=\eps^2\int\Ge(x,y)(-\Delta\varphi)(x)\dx+\int G(x,y)\varphi(x)\dx
\to\int\varphi\,\mathrm{d}\mu
\]
as $\eps\to0$. By density of $C^\infty$ in $C$, we conclude that
$\mu=\delta_y$. Proof of (ii). For $q=1$, we have:
\[
\|\Ge\ast h\|_1=\int|(\Ge\ast
h)(x)|\dx\le\int\dy|h(y)|\int\Ge(x,y)\dx=\int|h|=\|h\|_1.
\]
For $q=\infty$ we have, for any $x\in M$:
\[
|\Ge\ast
h(x)|\le\|h\|_\infty\int\Ge(x,y)\dy=\|h\|_\infty\int\Ge(x,y)\dx=\|h\|_\infty,
\]
and therefore $\|\Ge\ast h\|_\infty\le\|h\|_\infty$. The general
case follows by interpolation. Proof of (iii). Let $\Ue=\Ge\ast
h$. Then we can write
\[
-\eps^2\Delta(\Ue-h)+(\Ue-h)=\eps^2\Delta h.
\]
Multiplying by $|\Ue-h|^{q-2}(\Ue-h)$ and integrating, we obtain
\[
\eps^2(q-1)\int|\Ue-h|^{q-2}|\nabla(\Ue-h)|^2+\int|\Ue-h|^q=\eps^2\int\Delta
h|\Ue-h|^{q-2}(\Ue-h).
\]
By positivity of the first term above and H\"older's inequality,
\[
\int|\Ue-h|^q\le\eps^2\int|\Delta h||\Ue-h|^{q-1}\le\eps^2\|\Delta
h\|_q\|\Ue-h\|_q^{q-1}.
\]
Hence $\|\Ue-h\|_q\le\eps^2\|\Delta h\|_q$ and (iii) follows
recalling the definition of $\Ue$.
\end{proof}
Now we can prove Proposition~\ref{prop:subsol}.
\begin{proof}[Proof of Proposition~\ref{prop:subsol}]
Equation \eqref{fourthorder} is of the form:
\begin{equation}
\label{absfourth} \eps^2\Delta^2 u-\Delta u= \eps\lambda\,a(u)
+\lambda^2f'(e^{\sing+u})e^{\sing+u}(s-f(e^{\sing+u}))-\const,
\end{equation}
where $a$ is the operator defined by:
\begin{equation}
\label{adef}
a(u):=-[f''(e^{\sing+u})e^{\sing+u}+f'(e^{\sing+u})]e^{\sing+u}|\nabla(\sing+u)|^2
+2\Delta f(e^{\sing+u}).
\end{equation}
By \eqref{deltafid}, we can also write
\begin{align}
\label{adefsecond} a(u)=
[f''(e^{\sing+u})e^{\sing+u}+&f'(e^{\sing+u})]e^{\sing+u}|\nabla(\sing+u)|^2\\
\nonumber+&2f'(e^{\sing+u})e^{\sing+u}\Delta(\sing+u)
\end{align}
and therefore, recalling \eqref{singdelta} and \eqref{singnabla}
we can estimate:
\begin{equation}
\label{aest} \|a(u)\|_\infty\le\Phi(\|\Delta
u\|_\infty+\|u\|_{C^1}),
\end{equation}
for some continuous function $\Phi:[0,+\infty)\to\R$. We denote by
$\varphi$ a smooth function defined on $M$ with the following
properties:
\[
\varphi=
\begin{cases}
-A-1&\text{in}\ \cup_{j=1}^n B_\delta(p_j)\\
\varphi_0&\text{in}\ M\setminus\cup_{j=1}^n B_{2\delta}(p_j)\\
-A-1\le\varphi\le\varphi_0&\text{on}\ M\\
\int\varphi=0
\end{cases},
\]
where $\varphi_0$ is a suitable constant, and $\delta>0$ is
sufficiently small so that $B_{2\delta}(p_j)\cap
B_{2\delta}(p_k)=\emptyset$ for all $j,k=1,\ldots,n$ with $j\neq
k$. We denote by $\tue$ the unique solution for the problem
\begin{align}
\label{tuedef}
&\eps^2\Delta^2\tue-\Delta\tue=\varphi\qquad\text{on}\ M\\
\nonumber
&\int\tue=0.
\end{align}
Note that $\tue$ is well-defined. Indeed, since $\int\varphi=0$,
we have $\int\Ge\ast\varphi=0$. Therefore there exists a unique
solution $\tue$ for the problem $-\Delta\tue=\Ge\ast\varphi$
satisfying $\int\tue=0$. Writing $\eps^2\Delta^2-\Delta=
(-\eps^2\Delta+1)(-\Delta)$, we see that $\tue$ is the desired
unique solution for \eqref{tuedef}. By Lemma~\ref{lem:green}--(ii)
we have
$\|\Delta\tue\|_\infty=\|\Ge\ast\varphi\|_\infty\le\|\varphi\|_\infty$.
By elliptic regularity we have in turn
\begin{equation}
\label{c1}
\|\tue\|_{C^{1,\alpha}}\le C_1\|\varphi\|_\infty.
\end{equation}
We set $\uue:=\tue-k$,
where $k$ is defined by
\[
e^k=\frac{e^{\max_M\sing
+C_1\|\varphi\|_\infty}}{f^{-1}(\frac{s+f(0)}2)},
\]
and where $C_1$ is the constant in \eqref{c1}. In view of ($f$0)
and ($f$1), such a choice of $k$ implies that
$s-f(e^{\sing+\uue})\ge(s-f(0))/2>0$. Indeed, since $f$ is
strictly increasing, we have
\begin{align*}
f(e^{\sing+\uue})=f(e^{\sing+\tue-k})\le
f(e^{\max_M\sing+\|\tue\|_\infty-k})\le&
f(e^{\max_M\sing+C_1\|\varphi\|_\infty-k})\\
=&\frac{s+f(0)}{2}.
\end{align*}
 Now we check that for $\lambda\ge\lambda_0$ and
$\eps\le\eps_\lambda$, for suitable $\lambda_0$ and
$\eps_\lambda$, the function $\uue$ is indeed a subsolution for
\eqref{fourthorder}.
\medskip
\par
\textbf{ Claim:} There exists $\lambda_0>0$ such that for all
$\lambda\ge\lambda_0$ and for all $0<\eps<1$, $\uue$ is a
subsolution for \eqref{fourthorder} in $M\setminus\cup_{j=1}^n
B_\delta(p_j)$. Namely,
\begin{equation}
\label{subsolout}
 \varphi\le\eps\lambda\,a(\uue)+\lambda^2
f'(e^{\sing+\uue}) e^{\sing+\uue}(s-f(e^{\sing+\uue}))-A\quad
\text{in}\ M\setminus\cup_{j=1}^n B_\delta(p_j).
\end{equation}
Proof of \eqref{subsolout}. By \eqref{aest} and \eqref{c1} there
exists a constant $C_0$ such that
$\|a(\uue)\|_\infty+\|\uue\|_\infty+\max_M\sing\le C_0$. Let
$c_0=\min\{f'(t)\,/\,t\in[0,C_0]\}>0$ and
$\mu_0=\min_{M\setminus\cup_{j=1}^nB_\delta(p_j)}\sing$. It
suffices to check that:
\[
\varphi_0\le-\lambda C_0+\lambda^2c_0
e^{\mu_0-C_0}\frac{s-f(0)}{2} -A.
\]
The above inequality is clearly achieved for all
$\lambda\ge\lambda_0$, for sufficiently large $\lambda_0$. Hence,
\eqref{subsolout} is established.
\par
Now we fix $\lambda\ge\lambda_0$.
\medskip
\par
\textbf{ Claim:} For every fixed $\lambda\ge\lambda_0$, there
exists $\eps_\lambda>0$ such that $\uue$ is a subsolution for
\eqref{fourthorder} in $\cup_{j=1}^n B_\delta(p_j)$, for all
$0<\eps<\eps_\lambda$. Namely,
\begin{equation}
\label{subsolin}
\varphi\le\eps\lambda\,a(\uue)+\lambda^2f'(e^{\sing+\uue})e^{\sing+\uue}(s-f(e^{\sing+\uue}))
-A\quad \text{in}\ \cup_{j=1}^nB_\delta(p_j).
\end{equation}
Proof of \eqref{subsolin}. It suffices to prove the following
condition:
\[
-A-1\le-\eps\lambda C_0-A,
\]
which is clearly satisfied for all $0<\eps\le\eps_\lambda$, with
$\eps_\lambda>0$ such that $\eps_\lambda\lambda C_0\le1$. Hence,
\eqref{subsolin} is also established. Consequently, for
$\lambda\ge\lambda_0$ and for $0<\eps\le\eps_\lambda$, $\uue$ is a
subsolution for \eqref{fourthorder}, as asserted.
\par
We are left to analyze the asymptotic behavior of $\uue$ as
$\eps\to0$. We denote by $\tuo$ the unique solution for
$-\Delta\tuo=\varphi$ satisfying $\int\tuo=0$. Then $\tue-\tuo$
satisfies $-\Delta(\tue-\tuo)=\Ge\ast\varphi-\varphi$,
$\int(\tue-\tuo)=0$. By Lemma~\ref{lem:green}--(iii) we have
$\|\Ge\ast\varphi-\varphi\|_2\le\eps^2\|\Delta\varphi\|_2\to0$, as
$\eps\to0$. Therefore, $\|\tue-\tuo\|_{H^2}\le C\eps^2\to0$. It is
simple to check that $\uuo:=\tuo-k$ is a subsolution for
\eqref{gencs}. Clearly,
$\|\uue-\uuo\|_{H^2}=\|\tue-\tuo\|_{H^2}\to0.$
\end{proof}
Henceforth, $\lambda$ denotes a {\em fixed} constant satisfying
$\lambda\ge\lambda_0$.
\section{Existence of a local minimum}
\label{sec:locmin} We take $\lambda_0$ and $\eps_\lambda$ as in
Proposition~\ref{prop:subsol}. In this section we show:
\begin{Proposition}
\label{prop:locmin} Suppose $f$ satisfies ($f$0), ($f$1) and
($f$2). Then, (possibly taking a smaller $\eps_\lambda$), for
every fixed $\lambda\ge\lambda_0$ and for every
$0<\eps<\eps_\lambda$ there exists a solution $\ue$ for
\eqref{fourthorder}, corresponding to a local minimum for $\Ie$.
\end{Proposition}
We define
\[
\Ae:=\{u\in H^2\,/\ u\ge\uue\}.
\]
$\Ae$ is a closed convex subset of $H^2$ and therefore it is
weakly closed. It is readily checked that $\Ie$ attains its
minimum on $\Ae$, i.e., there exists $\ue$ such that
\[
\Ie(\ue)=\min_{\Ae}\Ie.
\]
The remaining part of this section is devoted to showing that
$\ue$ is a solution for \eqref{fourthorder} corresponding to a
local minimum for $\Ie$. The main issue is to show that $\ue$
belongs to the {\em interior} of $\Ae$ (in the sense of $H^2$),
and thus it is a critical point for $\Ie$. It is readily checked
that $\ue$ is a supersolution for \eqref{fourthorder}. Indeed, for
all $\phi\in H^2(M)$ such that $\phi\ge0$ and for all $t>0$ we
have $\ue+t\phi\in\Ae$, therefore
\[
\frac{\Ie(\ue+t\phi)-\Ie(\ue)}{t}\ge0.
\]
Consequently, taking into account \eqref{gateaux},
we obtain
\begin{align}
\label{uesupersol} 0\le\bra\Ie'(\ue),\phi\ket=
\eps^2&\int\Delta\ue\Delta\phi+\int\nabla\ue\cdot\nabla\phi\\
\nonumber
&+\eps\lambda\int\left[f''(e^{\sing+\ue})e^{\sing+\ue}+f'(e^{\sing+\ue})\right]
e^{\sing+\ue}|\nabla(\sing+\ue)|^2\phi\\
\nonumber
&+2\eps\lambda\int f'(e^{\sing+\ue})e^{\sing+\ue}\nabla(\sing+\ue)\cdot\nabla\phi\\
\nonumber &+\lambda^2\int
f'(e^{\sing+\ue})e^{\sing+\ue}(f(e^{\sing+\ue})-s)\phi
+\const\int\phi,
\end{align}
for all $\phi\in H^2$, $\phi\ge0$. Hence, $\ue$ is a supersolution
for \eqref{fourthorder}. We define
\[
\Ao=\{u\in H^1\,/u\ge\uuo\ \text{a.e.}\}.
\]
(Note that $\Ao$ is a subset of $H^1$, while $\Ae$ is a subset of
$H^2$). The next lemma provides estimates for $\ue$, independent
of $\eps\to0$. Throughout this section, we denote by $C>0$ a
general constant independent of $\eps$. Recall that $\Io$ is the
functional defined at the end of Section~\ref{sec:variational}.
\begin{Lemma}
\label{lem:uo} There exists a solution $\uo\in H^1$ for
\eqref{gencs} such that $\ue\to\uo$ strongly in $H^1$.
Furthermore,
\begin{itemize}
\item[(i)] $\lim_{\eps\to0}\Ie(\ue)=\inf_{\Ao}\Io$ \item[(ii)]
$\lim_{\eps\to0}\eps\|\Delta\ue\|_2=0$
\item[(iii)]$\lim_{\eps\to0}\eps\int
f'(e^{\sing+\ue})e^{\sing+\ue}|\nabla(\sing+\ue)|^2=0$.
\end{itemize}
\end{Lemma}
\begin{proof}
Since $\Ie(\ue)\le\Ie(\uue)\le C$ and since
$\int\ue\ge\int\uue\ge-C$, we readily have the following
estimates:
\begin{align}
\label{epsdeltabound}
&\eps\|\Delta\ue\|_2+\|\nabla\ue\|_2\le C\\
\nonumber
&\eps\int f'(e^{\sing+\ue})e^{\sing+\ue}|\nabla(\sing+\ue)|^2\le C\\
\nonumber &\left|\int\ue\right|\le C.
\end{align}
In particular, $\|\ue\|_{H^1}\le C$ and therefore we may assume that
for some $\uo\in H^1$ we have
$\ue\weak\uo$ weakly in $H^1$, strongly
in $L^p$ for all $1\le p<+\infty$ and a.e.
\par
Proof of (i). We can write
\[
\Ie(u)=\frac{\eps^2}{2}\|\Delta u\|_2^2
+\eps\lambda\int f'(e^{\sing+u})e^{\sing+u}|\nabla(\sing+u)|^2
+\Io(u)
\]
for all $u\in H^2$. In particular, $\Ie(u)\ge\Io(u)$ for all $u\in
H^2$. Since $\uue\to\uuo$ in $H^2$, we have
\begin{equation}
\label{infIo}
\inf_{\Ao}\Io=\inf_{\Ae}\Io+\circ_\eps(1),
\end{equation}
hence
\[
\Ie(\ue)=\inf_{\Ae}\Ie\ge\inf_{\Ae}\Io
=\inf_{\Ao}\Io+\circ_\eps(1).
\]
It follows that
\begin{equation}
\label{liminf}
\liminf_{\eps\to0}\Ie(\ue)\ge\inf_{\Ao}\Io.
\end{equation}
In order to obtain the inverse inequality, we fix $\eta>0$ and we
select $u_\eta\in\Ae$ such that
\[
\Io(u_\eta)\le\inf_{\Ae}\Io+\eta.
\]
(Note that $\Ae$ in not closed in $H^1$). We have:
\begin{align*}
\Ie(\ue)\le\Ie(u_\eta)\le&\Io(u_\eta)+\circ_\eps(1)\\
\le&\inf_{\Ae}\Io+\eta+\circ_\eps(1)\\
=&\inf_{\Ao}\Io+\eta+\circ_\eps(1).
\end{align*}
Therefore,
\begin{equation*}
\limsup_{\eps\to0}\Ie(\ue)\le\inf_{\Ao}\Io+\eta.
\end{equation*}
Since $\eta$ is arbitrary, we derive
\begin{equation}
\label{limsup}
\limsup_{\eps\to0}\Ie(\ue)\le\inf_{\Ao}\Io.
\end{equation}
Now by \eqref{liminf} and \eqref{limsup} the asserted equality (i)
is established.
\par
Proof of (ii) and (iii). By weak $H^1$ convergence and assumption
($f$2),
\[
\liminf_{\eps\to0}\Io(\ue)\ge\Io(\uo).
\]
Therefore, we have:
\begin{align*}
\inf_{\Ao}\Io=&\lim_{\eps\to0}\Ie(\ue)\\
=&\lim_{\eps\to0}\{\frac{\eps^2}{2}\|\Delta\ue\|_2^2
+\eps\lambda\int
f'(e^{\sing+\ue})e^{\sing+\ue}|\nabla(\sing+\ue)|^2
+\Io(\ue)\}\\
\ge&\lim_{\eps\to0}\{\frac{\eps^2}{2}\|\Delta\ue\|_2^2
+\eps\lambda\int f'(e^{\sing+\ue})e^{\sing+\ue}|\nabla(\sing+\ue)|^2\}+\Io(\uo)\\
\ge&\lim_{\eps\to0}\{\frac{\eps^2}{2}\|\Delta\ue\|_2^2
+\eps\lambda\int
f'(e^{\sing+\ue})e^{\sing+\ue}|\nabla(\sing+\ue)|^2\}+\inf_{\Ao}\Io.
\end{align*}
Therefore we obtain (ii) and (iii). Furthermore, we find that
$\lim_{\eps\to0}\Io(\ue)=\Io(\uo)$, which implies that $\ue\to\uo$
strongly in $H^1$, and that $\Io(\uo)=\inf_{\Ao}\Io$. Arguing
similarly as for $\ue$, we see that $\uo$ is a supersolution for
\eqref{gencs}. By the Hopf maximum principle,
\[
\uo>\uuo.
\]
Therefore, $\uo$ is a local minimum for $\Io$ in the $C^1$-topology.
By the Brezis-Nirenberg argument \cite{BN}, $\uo$ is a local minimum for $\Io$
in the $H^1$-topology and thus it is in fact a solution for \eqref{gencs}.
By elliptic regularity, $\uo$ is smooth.
\end{proof}
The next lemma shows that the strong maximum principle property
for $\uo$ and $\uuo$ carries over to $\ue$ and $\uue$, for small
values of $\eps$:
\begin{Lemma}
\label{lem:localmin} Suppose $f$ satisfies ($f$0), ($f$1) and
($f$2). Then, for all $\eps$ sufficiently small, $\ue>\uue$,
pointwise on $M$.
\end{Lemma}
\begin{proof}
We define
\[
\Fe=\eps\lambda\,a(\ue)
+\lambda^2f'(e^{\sing+\ue})e^{\sing+\ue}(s-f(e^{\sing+\ue}))-\const,
\]
where $a$ is the operator defined in \eqref{adef}. Then $\ue$
satisfies
\[
\eps^2\Delta^2\ue-\Delta\ue\ge\Fe.
\]
Since the Green function $\Ge$ for $-\eps^2\Delta+1$ is positive
(see Lemma~\ref{lem:green}), the above yields
\begin{equation}
\label{uesimple} -\Delta\ue\ge\Ge\ast\Fe.
\end{equation}
We define:
\[
\Fo=f'(e^{\sing+\uo})e^{\sing+\uo}(s-f(e^{\sing+\uo})) -\const.
\]
\medskip
\par
\textbf{Claim:} There exists some $q>1$ such
that
\begin{equation}
\label{rhsdecay} \|\Fe-\Fo\|_q\to0 \qquad\text{as}\ \eps\to0.
\end{equation}
Proof of \eqref{rhsdecay}. We show that $\eps a(\ue)\to0$. In view
of \eqref{adefsecond}, it suffices to show that
\[
\eps\|[f''(e^{\sing+\ue})e^{\sing+\ue}+f'(e^{\sing+\ue})]e^{\sing+\ue}|\nabla(\sing+\ue)|^2\|_q\to0
\]
and
\[
\eps\|f'(e^{\sing+\ue})e^{\sing+\ue}\Delta(\sing+\ue)\|_q\to0.
\]
By \eqref{welldef}, the fact $\|\ue\|_{H^1}\le C$ and Sobolev
embeddings we have
\[
\|[f''(e^{\sing+\ue})e^{\sing+\ue}+f'(e^{\sing+\ue})]e^{\sing+\ue}|\nabla(\sing+\ue)|^2\|_q\le
C
\]
for some $q>1$, and therefore the first limit follows easily. In
order to prove the second limit, we write, using
\eqref{singdelta}:
\[
\eps f'(e^{\sing+\ue})e^{\sing+\ue}\Delta(\sing+\ue) =-\eps
Af'(e^{\sing+\ue})e^{\sing+\ue}+\eps
f'(e^{\sing+\ue})e^{\sing+\ue}\Delta\ue.
\]
By similar arguments as above, the $L^q$-norm of the first term on
the right hand side above vanishes as $\eps\to0$. In order to
estimate the second term, we write for $r>2$ such that
$1/r+1/2=1/q$:
\begin{align*}
\|\eps
f'(e^{\sing+\ue})e^{\sing+\ue}\Delta\ue\|_q\le\|f'(e^{\sing+\ue})e^{\sing+\ue}\|_r\|\eps\Delta\ue\|_2
\le C\|\eps\Delta\ue\|_2\to0,
\end{align*}
where we used Lemma~\ref{lem:uo}--(ii) to derive the last step.
Hence, \eqref{rhsdecay} is established.
\par
By \eqref{rhsdecay} and by Lemma~\ref{lem:green}--(ii)--(iii), it
follows that
\begin{align*}
\|\Ge\ast\Fe-\Fo\|_q\le&\|\Ge\ast(\Fe-\Fo)\|_q+\|\Ge\ast\Fo-\Fo\|_q\\
\le&(\|\Fe-\Fo\|_q+\eps^2\|\Delta\Fo\|_q)\to0,
\end{align*}
as $\eps\to0$. We define $\we$ as the unique solution for
\[
(-\Delta+1)\we=\Ge\ast\Fe+\ue.
\]
Then, by \eqref{uesimple},
\[
(-\Delta+1)(\ue-\we)\ge0
\]
and therefore, by the maximum principle,
\[
\ue\ge\we.
\]
Since $\uo$ satisfies
\[
-\Delta\uo=\Fo
\]
we have
\[
(-\Delta+1)(\we-\uo)=\Ge\ast\Fe-\Fo+\ue-\uo
\]
and by standard elliptic estimates
\[
\|\we-\uo\|_{C^\alpha}\le C(\|\Ge\ast\Fe-\Fo\|_q+\|\ue-\uo\|_q)\to0.
\]
In conclusion, we have $\uue\to\uuo$ in $H^2$ and in particular
uniformly, $\we\to\uo$ uniformly and $\uo>\uuo$. It follows that
$\we>\uue$ for small $\eps$. Consequently, $\ue>\uue$ for small
$\eps$, as asserted.
\end{proof}
\begin{proof}[Proof of Proposition~\ref{prop:locmin}]
By Lemma~\eqref{lem:localmin} and the Sobolev embedding
$\|\phi\|_\infty\le C\|\phi\|_{H^2}$, $\ue$ belongs to the
interior of $\Ae$, in the sense of $H^2$. Therefore $\ue$ is a
critical point for $\Ie$ corresponding to a local minimum, as
asserted.
\end{proof}
\section{The Palais-Smale condition}
\label{sec:PS}
The main result of this section is:
\begin{Proposition}
\label{prop:ps} Suppose $f$ satisfies ($f$0), ($f$1) and ($f$3).
Then $\Ie$ satisfies the Palais-Smale condition.
\end{Proposition}
We denote by $(\uj)$, $\uj\in H^2$ a Palais-Smale sequence. That
is, we assume that $\Ie(\uj)\to\alpha$ for some $\alpha\in\R$ and
$\|\Ie'(\uj)\|_{H^{-1}}\to0$ as $j\to+\infty$. In order to prove
Proposition~\ref{prop:ps} we have to show that $(\uj)$ admits a
subsequence which converges strongly in $H^2$. By standard
compactness arguments, it suffices to show that $(\uj)$ is bounded
in $H^2$. It will be convenient to decompose
\[
\uj=\ujp+\cj,\qquad\qquad\int\ujp=0,\qquad\cj\in\R.
\]
Unless otherwise stated, throughout this section we denote by
$C>0$ a general constant independent of $j\to+\infty$, whose
actual value may vary from line to line. The Palais-Smale
assumption for $(\uj)$ implies in particular the following facts:
\begin{align}
\label{firstpsbound}
\Ie(\uj)=\frac{\eps^2}{2}&\int(\Delta\uj)^2+\frac{1}{2}\int|\nabla\uj|^2
+\eps\lambda\int f'(e^{\sing+\uj})e^{\sing+\uj}|\nabla(\sing+\uj)|^2\\
\nonumber
&+\frac{\lambda^2}{2}\int(f(e^{\sing+\uj})-s)^2+\const|M|\cj\to\alpha,
\end{align}
and (see \eqref{gateaux}):
\begin{align}
\label{test1}
\bra\Ie'(\uj),1\ket
=&\eps\lambda\int[f''(e^{\sing+\uj})e^{\sing+\uj}
+f'(e^{\sing+\uj})]e^{\sing+\uj}|\nabla(\sing+\uj)|^2\\
\nonumber
&+\lambda^2\int f'(e^{\sing+\uj})e^{\sing+\uj}(f(e^{\sing+\uj})-s)+A|M|
\le C,
\end{align}
\begin{align}
\label{testuj}
\bra\Ie'&(\uj),\uj\ket
=\eps^2\int(\Delta\uj)^2+\int|\nabla\uj|^2\\
\nonumber
&+\eps\lambda\int[f''(e^{\sing+\uj})e^{\sing+\uj}
+f'(e^{\sing+\uj})]e^{\sing+\uj}|\nabla(\sing+\uj)|^2\uj\\
\nonumber
&+2\eps\lambda\int f'(e^{\sing+\uj})e^{\sing+\uj}\nabla(\sing+\uj)\cdot\nabla\uj\\
\nonumber
&+\lambda^2\int f'(e^{\sing+\uj})e^{\sing+\uj}(f(e^{\sing+\uj})-s)\uj
+A|M|\cj\\
\nonumber
&\le\circ_j(1)(\|\Delta\uj\|_2+|\cj|),
\end{align}
\begin{align}
\label{testujp}
\bra\Ie'&(\uj),\ujp\ket
=\eps^2\int(\Delta\uj)^2+\int|\nabla\uj|^2\\
\nonumber
&+\eps\lambda\int[f''(e^{\sing+\uj})e^{\sing+\uj}
+f'(e^{\sing+\uj})]e^{\sing+\uj}|\nabla(\sing+\uj)|^2\ujp\\
\nonumber
&+2\eps\lambda\int f'(e^{\sing+\uj})e^{\sing+\uj}\nabla(\sing+\uj)\cdot\nabla\uj\\
\nonumber
&+\lambda^2\int f'(e^{\sing+\uj})e^{\sing+\uj}(f(e^{\sing+\uj})-s)\ujp
\le \circ_j(1)\|\Delta\uj\|_2.
\end{align}
\begin{Lemma}
\label{lem:psprelim} If either $\cj\ge-C$ or $\|\Delta\uj\|_2\le
C$, then $(\uj)$ is bounded in $H^2$.
\end{Lemma}
\begin{proof}
Suppose $\cj\ge-C$. Condition \eqref{firstpsbound} implies that
$\cj\le C$. Since $\cj\ge-C$, \eqref{firstpsbound} yields
$\|\Delta\uj\|_2\le C$. Hence $(\uj)$ is bounded in this case.
\par
Now suppose $\|\Delta\uj\|_2\le C$. Then, by Sobolev embeddings,
$\|\nabla\uj\|_p\le C_p$ for all $p\ge1$ and $\|\ujp\|_\infty\le
C$. It follows that
$\|e^\uj\|_\infty=e^{\cj}\|e^{\ujp}\|_\infty\le C$. Hence, by
\eqref{welldef}, we obtain $\int
f'(e^{\sing+\uj})e^{\sing+\uj}|\nabla(\sing+\uj)|^2\le C$.
Inserting into \eqref{firstpsbound}, we find that
$A|M|\cj+O_j(1)\to\alpha$. Consequently, $\cj\ge-C$, and therefore
$(\uj)$ is bounded also in this case.
\end{proof}
In view of Lemma~\ref{lem:psprelim}, we assume henceforth without
loss of generality:
\begin{align}
\label{ps1}
\|\Delta\uj\|_2\to+\infty
\qquad\text{as}\ j\to+\infty
\end{align}
and
\begin{equation}
\label{ps2}
\cj\to-\infty
\qquad\text{as}\ j\to+\infty.
\end{equation}
\begin{proof}
[Proof of Proposition~\ref{prop:ps} under assumption ($f$3)--(a)]
By ($f$1) there exists $t_0>0$ such that $f(t)-s>0$ for all $t\ge
t_0$, and consequently $f'(t)t(f(t)-s)\ge-C$ for some $C>0$
independent of $t$. We may further assume that $f'(t)t|f(t)-s|\le
f'(t)t(f(t)-s)+C$ and in view of ($f$3)--(a), we may assume that
$f^2(t)\le f'(t)t|f(t)|\le Cf'(t)t(f(t)-s)$ for all $t\ge t_0$.
Therefore, we derive from \eqref{test1} that
\begin{align}
\label{casei1}
&\int[f''(e^{\sing+\uj})e^{\sing+\uj}
+f'(e^{\sing+\uj})]e^{\sing+\uj}|\nabla(\sing+\uj)|^2\le C\\
\label{casei2}
&\int f'(e^{\sing+\uj})e^{\sing+\uj}|f(e^{\sing+\uj})-s|\le C\\
\label{casei3}
&\|f(e^{\sing+\uj})\|_2\le C.
\end{align}
Consequently, we may easily estimate the terms in the right hand
side of \eqref{testujp}. Indeed, using \eqref{casei1}, we have
\begin{align*}
\Big|\int[f''&(e^{\sing+\uj})e^{\sing+\uj}
+f'(e^{\sing+\uj})]e^{\sing+\uj}|\nabla(\sing+\uj)|^2\ujp\Big|\\
\le&\int[f''(e^{\sing+\uj})e^{\sing+\uj}
+f'(e^{\sing+\uj})]e^{\sing+\uj}|\nabla(\sing+\uj)|^2\,\|\ujp\|_\infty\\
\le&C\|\Delta\uj\|_2.
\end{align*}
Using \eqref{casei2}, we have
\begin{align*}
\Big|\int f'&(e^{\sing+\uj})e^{\sing+\uj}(f(e^{\sing+\uj})-s)\ujp\Big|\\
\le&\int f'(e^{\sing+\uj})e^{\sing+\uj}|f(e^{\sing+\uj})-s|\,\|\ujp\|_\infty\\
\le&C\|\Delta\uj\|_2.
\end{align*}
Finally, integrating by parts, we have:
\begin{align*}
\int f'&(e^{\sing+\uj})e^{\sing+\uj}\nabla(\sing+\uj)\cdot\nabla\uj\\
=&\int\nabla f(e^{\sing+\uj})\cdot\nabla\uj\\
=&-\int f(e^{\sing+\uj})\Delta\uj.
\end{align*}
Hence, using \eqref{casei3} we have
\begin{align*}
\Big|\int f'&(e^{\sing+\uj})e^{\sing+\uj}\nabla(\sing+\uj)\cdot\nabla\uj\Big|\\
=&\Big|\int f(e^{\sing+\uj})\Delta\uj\Big|\\
\le&\|f(e^{\sing+\uj})\|_2\|\Delta\uj\|_2
\le C\|\Delta\uj\|_2.
\end{align*}
Inserting into \eqref{testujp} we find $\|\Delta\uj\|_2\le C$,
which is in contradiction with \eqref{ps1}.
\end{proof}
In order to prove Proposition~\ref{prop:ps} in the remaining case
($f$3)--(b), we first establish an identity:
\begin{Lemma}
\label{lem:identity} For all $u\in H^2$ the following identity
holds:
\begin{align*}
\int[f''&(e^{\sing+u})e^{\sing+u}
+f'(e^{\sing+u})]e^{\sing+u}|\nabla(\sing+u)|^2u\\
&+2\int f'(e^{\sing+u})e^{\sing+u}\nabla(\sing+u)\cdot\nabla u\\
=&\int f'(e^{\sing+u})e^{\sing+u}\nabla(\sing+u)\cdot\nabla u-\int
f'(e^{\sing+u})e^{\sing+u}\Delta(\sing+u)u.
\end{align*}
\end{Lemma}
\begin{proof}
Integrating by parts, we have:
\begin{align*}
\int[f''&(e^{\sing+u})e^{\sing+u}
+f'(e^{\sing+u})]e^{\sing+u}|\nabla(\sing+u)|^2u\\
=&\int\nabla[f'(e^{\sing+u})e^{\sing+u}]\cdot\nabla(\sing+u)u\\
=&-\int f'(e^{\sing+u})e^{\sing+u}\Delta(\sing+u)u -\int
f'(e^{\sing+u})e^{\sing+u}\nabla(\sing+u)\cdot\nabla u.
\end{align*}
The asserted identity follows.
\end{proof}
\begin{proof}[Proof of Proposition~\ref{prop:ps} under assumption ($f$3)--(b)]
By Lemma~\ref{lem:identity} with $u=\uj$, condition \eqref{testuj}
may be equivalently written in the form:
\begin{align*}
\bra\Ie'&(\uj),\uj\ket
=\eps^2\int(\Delta\uj)^2+\int|\nabla\uj|^2\\
&+\eps\lambda\int f'(e^{\sing+\uj})e^{\sing+\uj}\nabla(\sing+\uj)\cdot\nabla\uj\\
&-\eps\lambda\int f'(e^{\sing+\uj})e^{\sing+\uj}\Delta(\sing+\uj)\uj\\
&+\lambda^2\int f'(e^{\sing+\uj})e^{\sing+\uj}(f(e^{\sing+\uj})-s)\uj
+A|M|\cj\\\nonumber
&\le\circ_j(1)(\|\Delta\uj\|_2+|\cj|),
\end{align*}
and consequently, we have
\begin{align}
\label{testujequiv}
\eps^2\int(\Delta\uj)^2+\int|\nabla\uj|^2
&\le -\eps\lambda\int f'(e^{\sing+\uj})e^{\sing+\uj}\nabla(\sing+\uj)\cdot\nabla\uj\\
\nonumber
&+\eps\lambda\int f'(e^{\sing+\uj})e^{\sing+\uj}\Delta(\sing+\uj)\uj\\
\nonumber
&-\lambda^2\int f'(e^{\sing+\uj})e^{\sing+\uj}(f(e^{\sing+\uj})-s)\uj
-A|M|\cj\\
\nonumber
&+\circ_j(1)(\|\Delta\uj\|_2+|\cj|).
\end{align}
We estimate term-by-term the right hand side in
\eqref{testujequiv}.
\medskip
\par
\textbf{Claim:} There holds:
\begin{equation}
\label{rhsest1}
-\int f'(e^{\sing+\uj})e^{\sing+\uj}\nabla(\sing+\uj)\cdot\nabla\uj\le C\|\Delta\uj\|_2.
\end{equation}
Proof of \eqref{rhsest1}. Since $f'(t)t\le C$ we have, for any
$1<p<2$:
\begin{align*}
&-\int f'(e^{\sing+\uj})e^{\sing+\uj}\nabla(\sing+\uj)\cdot\nabla\uj\\
=&-\int f'(e^{\sing+\uj})e^{\sing+\uj}\nabla\sing\cdot\nabla\uj
-\int f'(e^{\sing+\uj})e^{\sing+\uj}|\nabla\uj|^2\\
\le&|\int f'(e^{\sing+\uj})e^{\sing+\uj}\nabla\sing\cdot\nabla\uj|\\
\le&C\|\nabla\sing\|_p\|\nabla\uj\|_{p'}
\le C\|\Delta\uj\|_2,
\end{align*}
and thus \eqref{rhsest1} is established.
\medskip
\par
\textbf{Claim:} There holds:
\begin{equation}
\label{rhsest2} |\int
f'(e^{\sing+\uj})e^{\sing+\uj}\Delta(\sing+\uj)\uj| \le
C\|\Delta\uj\|_2.
\end{equation}
Proof of \eqref{rhsest2}.
By the assumption $\sup_{t>0}f'(t)t|\log t|\le+\infty$,
we readily derive:
\begin{align*}
|\int&f'(e^{\sing+\uj})e^{\sing+\uj}\Delta\uj\uj|\\
\le&|\int f'(e^{\sing+\uj})e^{\sing+\uj}(\sing+\uj)\Delta\uj|
+|\int f'(e^{\sing+\uj})e^{\sing+\uj}\sing\Delta\uj|\\
\le&C\|\Delta\uj\|_2+C\|\sing\|_2\|\Delta\uj\|_2\\
\le&C\|\Delta\uj\|_2.
\end{align*}
On the other hand, recalling \eqref{singdelta} we have:
\begin{align*}
\int f'(e^{\sing+\uj})&e^{\sing+\uj}\Delta\sing\,\uj =-A\int
f'(e^{\sing+\uj})e^{\sing+\uj}\uj \\=&-A\int
f'(e^{\sing+\uj})e^{\sing+\uj}(\sing+\uj)+A\int
f'(e^{\sing+\uj})e^{\sing+\uj}\sing.
\end{align*}
Hence, in view of ($f$3)--(b) we derive:
\begin{align*}
|\int f'(e^{\sing+\uj})e^{\sing+\uj}\Delta\sing\,\uj|\le C.
\end{align*}
Now \eqref{rhsest2} follows (recall that by assumption
$\|\Delta\uj\|_2\to+\infty$).
\par
Finally, we easily estimate using ($f$3)--(b):
\begin{align*}
|\int f'(e^{\sing+\uj})e^{\sing+\uj}(f(e^{\sing+\uj}-s))\uj| \le
C\int|\uj|\le C(\|\Delta\uj\|_2+|\cj|).
\end{align*}
Inserting \eqref{rhsest1}--\eqref{rhsest2} and the above estimate
into \eqref{testujequiv}, we finally obtain
\begin{equation}
\label{testujconcl} \eps^2\|\Delta\uj\|_2^2\le
C(\|\Delta\uj\|_2+|\cj|).
\end{equation}
In turn, \eqref{testujconcl} yields
\begin{equation}
\label{normprelim}
\|\Delta\uj\|_2\le C|\cj|^{1/2}.
\end{equation}
Using \eqref{normprelim}, we show:
\bigskip
\par
\textbf{Claim:} There holds $\|e^\uj\|_\infty\to0$. More
precisely, the following estimate holds, for some $\gamma>0$:
\begin{equation}
\label{euest}
\|e^\uj\|_\infty\le e^{-\gamma|\cj|}.
\end{equation}
Proof of \eqref{euest}.
By \eqref{normprelim} and Sobolev embeddings, we have:
\begin{align*}
\|e^\uj\|_\infty\le e^{\cj+\|\ujp\|_\infty}\le e^{\cj+C\|\Delta\uj\|_2}
\le e^{-|\cj|+C|\cj|^{1/2}},
\end{align*}
and \eqref{euest} follows. \par Finally, we show:
\bigskip
\par
\textbf{Claim:} There holds:
\begin{equation}
\label{badtermest}
|\int f'(e^{\sing+\uj})e^{\sing+\uj}\Delta(\sing+\uj)\ujp|
=\circ_j(1)\|\Delta\uj\|_2.
\end{equation}
Proof of \eqref{badtermest}. By \eqref{euest} and the fact
$\|\ujp\|_\infty\le C\|\Delta\uj\|_2\le C|\cj|^{1/2}$, we have
\begin{align*}
|\int f'(e^{\sing+\uj})e^{\sing+\uj}\Delta\uj\ujp|
\le&\|f'(e^{\sing+\uj})e^{\sing+\uj}\|_\infty\|\ujp\|_\infty\|\Delta\uj\|_2\\
\le&C\|e^\uj\|_\infty\|\ujp\|_\infty\|\Delta\uj\|_2\\
\le&Ce^{-\gamma|\cj|}|\cj|^{1/2}\|\Delta\uj\|_2\\
=&\circ_j(1)\|\Delta\uj\|_2.
\end{align*}
On the other hand, recalling \eqref{singdelta}, we have
by similar arguments as above:
\begin{align*}
|\int f'(e^{\sing+\uj})e^{\sing+\uj}\Delta\sing\ujp|
=&A|\int f'(e^{\sing+\uj})e^{\sing+\uj}\ujp|\\
\le&C e^{-\gamma|\cj|}|\cj|^{1/2}\to0.
\end{align*}
Hence \eqref{badtermest} is established.
\par
Now, using Lemma~\ref{lem:identity} with $u=\ujp$, we rewrite
\eqref{testujp} in the form:
\begin{align*}
\eps^2\int(\Delta\uj)^2+&\int|\nabla\uj|^2
\le-\eps\lambda\int f'(e^{\sing+\uj})e^{\sing+\uj}\nabla(\sing+\uj)\cdot\nabla\uj\\
&+\eps\lambda\int f'(e^{\sing+\uj})e^{\sing+\uj}\Delta(\sing+\uj)\ujp\\
&-\lambda^2\int
f'(e^{\sing+\uj})e^{\sing+\uj}(f(e^{\sing+\uj})-s)\ujp
+\circ_j(1)\|\Delta\uj\|_2,
\end{align*}
In view of \eqref{rhsest1}, \eqref{badtermest} and assumption
[$f$3]-b, we derive from the above that $\|\Delta\uj\|_2\le C$.
This is in contradiction with \eqref{ps1}.
\end{proof}
Now we can finally prove our main result:
\begin{proof}[Proof of Theorem~\ref{thm:main}]
By Proposition~\ref{prop:locmin}, the functional $\Ie$ admits a
critical point corresponding to a local minimum. By
Proposition~\ref{prop:ps}, $\Ie$ satisfies the Palais-Smale
condition. If $\ue$ is not a strict local minimum, it is known
that $\Ie$ has a continuum of critical points (see, e.g.,
\cite{T}). In particular, $\Ie$ has at least two critical points.
If $\ue$ is a strict local minimum, we note that for $c\in\R$,
$c\to-\infty$, we have $\Ie(c)\to-\infty$. Therefore $\Ie$ admits
a mountain pass structure in the sense of Ambrosetti and
Rabinowitz \cite{AR}. Hence by the ``mountain pass theorem"
\cite{AR} we obtain the existence of a second critical point for
$\Ie$. In either case, we conclude that the fourth order equation
\eqref{fourthorder}, (equivalently, system
\eqref{genmcsa'}--\eqref{genmcsb'}) admits at least two solutions,
as asserted.
\end{proof}
\section*{Acknowledgments}
In am grateful to Professor Gabriella Tarantello for interesting
and stimulating discussions.

\end{document}